\documentclass{article}
\usepackage{amssymb,latexsym,amsmath,amsthm,mathrsfs}
\usepackage{graphicx}
\newcommand{\comment}[1]{}

\begin{document}
\title{Finding the sum of any series from a given general term\footnote{Presented
to the St. Petersburg Academy on October 13, 1735.
Originally
published as
{\em Inventio summae cuiusque seriei ex dato termino generali},
Commentarii academiae scientiarum Petropolitanae \textbf{8} (1741), 9--22.
E47 in the Enestr{\"o}m index.
Translated from the Latin by Jordan Bell,
Department of Mathematics, University of Toronto, Toronto, Ontario, Canada.
Email: jordan.bell@gmail.com}}
\author{Leonhard Euler}
\date{}

\maketitle

1. When I had considered more carefully what I explained by the
geometrical method in the previous paper\footnote{Translator: {\em Methodus universalis serierum convergentium summas quam proxime inveniendi}, E46.}
on the summation of series and when I had investigated the same method of summation
analytically, I saw that what I had extracted geometrically could
be deduced from a special method of summation that I had already mentioned
three years earlier in a paper on the summation of series\footnote{Translator: 
{\em Methodus generalis summandi progressiones}, E25, \S 2.}.
But I had not thought about this more since then.
Having examined more deeply the effectiveness of the analytical method,
I perceived that not only was the formula discovered geometrically
contained in it, but also that by means of it more could be accomplished by
adding more terms, so that it would show the true sum absolutely.
The geometrical method however seems to find these same terms with the greatest
difficulty.

2. In the former paper on the summation of series, if $x$ is the general term
of index $n$ of some series, I exhibited in a universal way the following form
for the summatory term
\[
\int x dn+\frac{x}{2}+\frac{dx}{12dn}-\frac{d^3x}{720dn^3}+\textrm{etc.},
\]
in which the differentials of $x$ over powers of the differential $dn$, which
is assumed constant, are destroyed, because $x$ is taken to be given by
$n$,\footnote{Translator: My best reading is that since $dn$ is small but
fixed, if $d^kx=0$ for some $k$ then $\frac{d^k x}{dn^k}=0$ and also for all higher
powers.} 
so that an algebraic sum is obtained if of course $xdn$ admits integration.
In the integration of $xdn$ indeed a constant ought to be added such that the whole
expression vanishes by putting $n=0$.

3. Now since I have set out in this paper to describe more accurately 
this formula and its use, before everything else I shall explain how I
discovered the formula. I used some singular arguments which offer
much to Analysis,  
partly new and partly already known, which however as far as I recall
are not demonstrated clearly enough elsewhere.

4. It follows from the nature of infinitesimal calculus that if $y$ depends
in any fixed way on $x$, if $x+dx$ is put in place of $x$ then $y$ will turn into
$y+dy$. Now, if $x$ is then increased by the element $dx$, 
that is $x$ is changed to $x+2dx$,
then in place of $y$ we will have $y+2dy+ddy$.
And if $x$ is again increased with $dx$,
then $y$ will transform into $y+3dy+3ddy+d^3y$, where
the coefficients are the same as those of the powers of a
binomial.
From here
it follows that if $x+mdx$ is put in place of $x$, then $y$ will take on this form
\[
y+\frac{m}{1}dy+\frac{m(m-1)}{1\cdot 2}ddy+\frac{m(m-1)(m-2)}{1\cdot 2\cdot 3}d^3 y+\textrm{etc.}
\]

5. Now for our purpose let $m$ be an infinitely large number such that
$mdx$ represents a finite quantity; putting $x+mdx$ in place of $x$, $y$ will
have this value
\[
y+\frac{mdy}{1}+\frac{m^2 d^2y}{1\cdot 2}+\frac{m^3d^3y}{1\cdot 2\cdot 3}+
\frac{m^4d^4y}{1\cdot 2\cdot 3\cdot 4}+\textrm{etc.}
\]
Now if we let $mdx=a$ or $m=\frac{a}{dx}$, if $x+a$ is put for $x$, then
$y$ will
assume this form
\[
y+\frac{ady}{1dx}+\frac{a^2ddy}{1\cdot 2dx^2}+\frac{a^3d^3y}{1\cdot 2\cdot 3dx^3}+\textrm{etc.},
\]
in which all the terms are of finite magnitude.

6. This series, which exhibits the transformed value of $y$ if $x+a$ is put
in place of $x$, was first found by the very insightful Taylor in
the {\em Methodus incrementorum directa et inversa}, and he applied it
to many excellent uses. The first result that follows is the raising of a binomial to any power. So if the value of $(x+a)^m$ is sought, I put
\[
y=x^m
\]
and if $x+a$ is put in place of $x$, the value of $y$ will be $(x+a)^m$.
Since therefore
\[
dy=mx^{m-1}dx, \quad d^2y=m(m-1)x^{m-2}dx^2
\]
and so on, it will be
\[
(x+a)^m=x^m+\frac{max^{m-1}}{1}+\frac{m(m-1)a^2x^{m-2}}{1\cdot 2}+\textrm{etc.}
\]

7. Then by doing the following the Taylor series lets us find approximately
a root of this equation. Let us have an equation involving an unknown
$z$, namely $Z=0$, where $Z$ is a quantity composed in some known way
from the unknown $z$.
Then take $x$ as a value nearly equal to $z$, and let the quantity
of $Z$ which occurs when $x$ is put in place of $z$ be put $=y$,
so that if $x$ were the true value of $z$ then $y=0$.

8. Now since $x$ differs from the true value of $z$ by a certain amount,
put the true value of $z$ to be $x+a$. It is thus clear that if in $y$ we put
$x+a$ in place of $x$ then $y$ will vanish. And indeed if one puts $x+a$
in place of $x$ then $y$ will turn into
\[
y+\frac{ady}{1dx}+\frac{a^2ddy}{1\cdot 2dx^2}+\frac{a^3d^3y}{1\cdot 2\cdot 3dx^3}+\textrm{etc.}
\]
From this it follows that
\[
0=y+\frac{ady}{1dx}+\frac{a^2ddy}{1\cdot 2dx^2}+\textrm{etc.}
\]

9. Since $x$ is set to be very close to $z$, $a$ will be a very small quantity,
so that beside the first two terms all the following ones will vanish.
By doing this it arises that $a=-\frac{ydx}{dy}$ and so
$z=x-\frac{ydx}{dy}$,
which is a value 
much nearer to $z$ than $x$. 
Thus for the equation
\[
z^3-3z-20=0
\]
it will be
\[
y=x^3-3x-20 \quad \textrm{and} \quad \frac{dy}{dx}=3x^2-3
\]
and hence
\[
z=x-\frac{x^3-3x-20}{3xx-3}=\frac{2x^3+20}{3xx-3}.
\]
Now by first taking $x=3$ it will be $z=3\frac{1}{12}$, and then repeating
this for a second time taking this value in place of $x$ will lead
to a value even closer to $z$.

10. Next, if some condition is stipulated on the function $y$ by which it is to have a particular relation to $x$, then the above formula will turn into
an equation which contains the character of $y$. Thus if $y$ is a function of
$x$ that vanishes by putting $x=0$, I put $a=-x$;
for thus it turns out that $x+a=0$ and it will be
\[
0=y-\frac{xdy}{1dx}+\frac{x^2ddy}{1\cdot 2dx^2}-\frac{x^3d^3y}{1\cdot 2\cdot 3dx^3}+\textrm{etc.}
\]
or
\[
y=\frac{xdy}{1dx}-\frac{x^2ddy}{1\cdot 2dx^2}+\frac{x^3d^3y}{1\cdot 2\cdot 3dx^3}-\textrm{etc.}
\]
The nature of all functions of $x$ which vanish by putting $x=0$ are contained in this equation.

11. If we write $\int zdx$ for $y$, it will be
\[
dy=zdx, \quad ddy=dzdx, \quad d^3y=d^2zdx \quad \textrm{etc.};
\]
substituting these values in we get
\[
\int zdx=\frac{xz}{1}-\frac{x^2dz}{1\cdot 2dx}+\frac{x^3ddz}{1\cdot 2\cdot 3dx^2}-\textrm{etc.},
\]
in which equation the integral of $zdx$ is expressed by an infinite series.
And this is the general quadrature of curves which the most insightful Johann Bernoulli
gave in the {\em Acta eruditorum} of Leipzig;
however, he did not attach the analysis which led to this series.

12. However disregarding this, which pertains less to our purpose,
I return to series. Therefore let us have some series
\[
A+B+C+D+\cdots+X,
\]
in which $A$ denotes the first term, $B$ the second, and $X$ that whose
index is $x$, so that $X$ is the general term of the given series.
Let us also put the sum of this progression to be
\[
A+B+C+D+\cdots+X=S;
\]
$S$ will be the summatory term, and if the series is determined it
will be composed from $x$ and fixed as much as $X$ is.

13. Now because $S$ will exhibit the sum of as many terms from the series
as there are unities in $x$, if $x-1$ is written in place of $x$ in $S$,
we will obtain the previous sum with the final term $X$ removed.
This substitution therefore turns $S$ into $S-X$.
Let us compare this with the above formula;
it will be $S=y$ and $a=-1$, from which the transformed value of $S$, or $S-X$, it will be
\[
=S-\frac{dS}{1dx}+\frac{ddS}{1\cdot 2dx^2}-\frac{d^3S}{1\cdot 2\cdot 3dx^3}+\textrm{etc.},
\]
from which this equation arises:
\[
X=\frac{dS}{1dx}-\frac{ddS}{1\cdot 2dx^2}+\frac{d^3S}{1\cdot 2\cdot 3dx^3}
-\frac{d^4S}{1\cdot 2\cdot 3\cdot 4dx^4}+\textrm{etc.}
\]

14. Therefore by means of this equation, the general term of any series
is found from the given summatory term.
However, since this is already very easy, it would be superfluous
to use this method for finding the general term from the summatory term.
Rather this equation is most useful if all the terms are expanded,
and it can thus be applied to all uses.
For by a known method the series
\[
X=\frac{dS}{1dx}-\frac{ddS}{1\cdot 2dx^2}+\frac{d^3S}{1\cdot 2\cdot 3dx^3}
-\textrm{etc.}
\]
can be inverted, so that from the general term $X$ the summatory
term $S$ can be determined, which is desired most.

15. Let us therefore put
\[
\frac{dS}{dx}=\alpha X+\frac{\beta dX}{dx}+\frac{\gamma ddX}{dx^2}
+\frac{\delta d^3X}{dx^3}+\frac{\epsilon d^4X}{dx^4}+\textrm{etc.},
\]
so it will be
\[
S=\alpha \int Xdx+\beta X+\frac{\gamma dX}{dx}+\frac{\delta ddX}{dx^2}+\textrm{etc.}
\]
Next, it will be
\[
\frac{ddS}{dx^2}=\frac{\alpha dX}{dx}+\frac{\beta ddX}{dx^2}
+\frac{\gamma d^3X}{dx^3}+\frac{\delta d^4X}{dx^4}+\textrm{etc.}
\]
and
\[
\frac{d^3S}{dx^3}=\frac{\alpha ddX}{dx^2}+\frac{\beta d^3X}{dx^3}+\frac{\gamma d^4X}{dx^4}+\textrm{etc.}
\]
and
\[
\frac{d^4S}{dx^4}=\frac{\alpha d^3X}{dx^3}+\frac{\beta d^4X}{dx^4}+\textrm{etc.}
\]
and then
\[
\frac{d^5S}{dx^5}=\frac{\alpha d^4X}{dx^4}+\textrm{etc.}
\]

16. Let us substitute the series
in place of each of the terms of the above series, and put the similar
terms among these equal to $0$. By doing this, the coefficients
$\alpha,\beta,\gamma$ etc. will be determined as follows

\begin{eqnarray*}
\alpha&=&1,\\
\beta&=&\frac{\alpha}{2},\\
\gamma&=&\frac{\beta}{2}-\frac{\alpha}{6},\\
\delta&=&\frac{\gamma}{2}-\frac{\beta}{6}+\frac{\alpha}{24},\\
\epsilon&=&\frac{\delta}{2}-\frac{\gamma}{6}+\frac{\beta}{24}-\frac{\alpha}{120},\\
\zeta&=&\frac{\epsilon}{2}-\frac{\delta}{6}+\frac{\gamma}{24}-\frac{\beta}{120}+\frac{\alpha}{720}\\
&&\textrm{etc.}
\end{eqnarray*}

17. Thus the coefficients $\alpha,\beta,\gamma,\delta$ etc. constitute a series
of such a nature that each term is determined by all the preceding terms, with the first term being $=1$. Also, the numbers which all the final terms 
need to be divided by constitute the progression called {\em hypergeometric} by Wallis
\[
2, \quad 6, \quad 24, \quad 120, \quad 720, \quad 5040 \quad \textrm{etc.}
\]
However, this series of coefficients $\alpha,\beta,\gamma$ etc. is thus
constituted that I could hardly believe that each could be exhibited by some
general term.

18. Therefore for our purposes we should be contented with the series
of coefficients being 
continued as far as we want, which can easily be done perfectly
from the law of the progression. I have worked out this series as follows,
\[
\begin{split}
&+1, \, +\frac{1}{1\cdot 2}, \, +\frac{1}{1\cdot 2\cdot 3\cdot 2}, \, +0, \,
-\frac{1}{1\cdot 2\cdot 3\cdot 4\cdot 5\cdot 6}, \, -0,
\,+\frac{1}{1\cdot 2\cdot 3\cdot 4\cdot 5\cdot 6\cdot 7\cdot 6}, \,+0,\\
&-\frac{3}{1\cdot 2\cdot 3\cdot 4\cdot 5\cdot 6\cdot 7\cdot 8\cdot 9\cdot 10},
\, -0, \, +\frac{5}{1\cdots 11\cdot 6}, \, +0, \, -\frac{691}{1\cdots 13\cdot 210},
\, -0,\\
&+\frac{35}{1\cdots 15\cdot 2}, +0, -\frac{3617}{1\cdots 17\cdot 30} \quad
\textrm{etc.}
\end{split}
\] 
It is notable that in this series all the even terms besides the second vanish.

19. Therefore if these terms are substituted in place of
$\alpha,\beta,\gamma$ etc., we will obtain the following summatory term
\begin{eqnarray*}
S&=&\int Xdx+\frac{X}{1\cdot 2}+\frac{dX}{1\cdot 2\cdot 3\cdot 2dx}-\frac{d^3X}{1\cdot 2\cdot 3\cdot 4\cdot 5\cdot 6dx^3}
+\frac{d^5X}{1\cdots 7\cdot 6dx^5}\\
&&-\frac{3d^7X}{1\cdots 9\cdot 10dx^7}+\frac{5d^9X}{1\cdots 11\cdot 6dx^9}
-\frac{691d^{11}X}{1\cdots 13\cdot 210dx^{11}}\\
&&+\frac{35d^{13}X}{1\cdots 15\cdot 2dx^{13}}-\frac{3617d^{15}X}{1\cdots 17\cdot 30dx^{15}}+\textrm{etc.}
\end{eqnarray*}

20. This series has an important use in finding the sums of algebraic progressions, in which $x$ does not appear in the denominator of the general term.
For  this reason $x$ will have positive exponents everywhere, and hence
some differential of it will vanish and thus the series will stop, and therefore
the summatory term will be represented by a finite number of terms. 
Immediately we see that all the terms which do not contain $x$ can be ignored,
since already some constant needs to be added in $\int Xdx$, to make $S=0$
when we put $x=0$.

21. To clearly see the use of this formula, it is worthwhile to offer
some examples. Thus let $X=x$, that is, let the series to be summed be
\[
1+2+3+\cdots+x;
\]
since
\[
\int X dx=\frac{x^2}{2}
\]
the sum will be
\[
S=\frac{x^2+x}{2};
\]
for $\frac{dX}{dx}$ is constant and is therefore ignored, and the following
differentials spontaneously vanish.

Next let $X=x^2$, or let this be the series to be summed
\[
1+4+9+\cdots+x^2;
\]
it will be
\[
\int Xdx=\frac{x^3}{3} \quad \textrm{and} \quad \frac{dX}{dx}=2x
\]
and hence the sum of the series is
\[
S=\frac{x^3}{3}+\frac{x^2}{2}+\frac{x}{6}.
\]

22. Now let the general series of powers of the natural numbers be given
\[
1+2^n+3^n+4^n+5^n+\textrm{etc.},
\]
whose general term is $x^n$. One will therefore have $X=x^n$ and
\[
\int Xdx=\frac{x^{n+1}}{n+1}.
\]
Furthermore the differentials will be thus obtained,
\begin{eqnarray*}
\frac{dX}{dx}&=&nx^{n-1},\\
\frac{d^3X}{dx^3}&=&n(n-1)(n-2)x^{n-3},\\
\frac{d^5X}{dx^5}&=&n(n-1)(n-2)(n-3)(n-4)x^{n-5}\\
&&\textrm{etc.}
\end{eqnarray*}
Therefore with these values substituted the summatory term of the given series
will be
\begin{eqnarray*}
S&=&\frac{x^{n+1}}{n+1}+\frac{x^n}{2}+\frac{nx^{n-1}}{2\cdot 6}
-\frac{n(n-1)(n-2)x^{n-3}}{2\cdot 3\cdot 4\cdot 30}
+\frac{n(n-1)(n-2)(n-3)(n-4)x^{n-5}}{2\cdot 3\cdot 4\cdot 5\cdot 6\cdot 42}\\
&&-\frac{n(n-1)\cdots(n-6)x^{n-7}}{2\cdot 3\cdots 8\cdot 30}
+\frac{n(n-1)\cdots(n-8)5x^{n-9}}{2\cdot 3\cdots 10\cdot 66}
-\frac{n(n-1)\cdots(n-10)691x^{n-11}}{2\cdot 3\cdots 12\cdot 2730}\\
&&+\frac{n(n-1)\cdots(n-12)7x^{n-13}}{2\cdot 3\cdots 14\cdot 6}
-\frac{n(n-1)\cdots(n-14)3617 x^{n-15}}{2\cdot 3\cdots 16\cdot 510}
+\textrm{etc.}
\end{eqnarray*}
The above series $\alpha,\beta,\gamma$ etc. should be continued as far
necessary for this series, which is worth continuing.

23. Thus from this general summation of the series whose general term is $x^n$,
sums of series of particular powers can be constructed, as follows,
\begin{eqnarray*}
\int x^1&=&\frac{x^2}{2}+\frac{x}{2},\\
\int x^2&=&\frac{x^3}{3}+\frac{x^2}{2}+\frac{x}{6},\\
\int x^3&=&\frac{x^4}{4}+\frac{x^3}{2}+\frac{x^2}{4},\\
\int x^4&=&\frac{x^5}{5}+\frac{x^4}{2}+\frac{x^3}{3}-\frac{x}{30},\\
\int x^5&=&\frac{x^6}{6}+\frac{x^5}{2}+\frac{5x^4}{12}-\frac{x^2}{12},\\
\int x^6&=&\frac{x^7}{7}+\frac{x^6}{2}+\frac{x^5}{2}-\frac{x^3}{6}+\frac{x}{42},\\
\int x^7&=&\frac{x^8}{8}+\frac{x^7}{7}+\frac{7x^6}{12}-\frac{7x^4}{24}+\frac{x^2}{12},\\
\int x^8&=&\frac{x^9}{9}+\frac{x^8}{2}+\frac{2x^7}{3}-\frac{7x^5}{15}+\frac{2x^3}{9}-\frac{x}{30},\\
\int x^9&=&\frac{x^{10}}{10}+\frac{x^9}{2}+\frac{3x^8}{4}-\frac{7x^6}{10}+\frac{x^4}{2}-\frac{3x^2}{20},\\
\int x^{10}&=&\frac{x^{11}}{11}+\frac{x^{10}}{2}+\frac{5x^9}{6}-x^7+x^5-\frac{x^3}{2}+\frac{5x}{66},\\
\int x^{11}&=&\frac{x^{12}}{12}+\frac{x^{11}}{2}+\frac{11x^{10}}{12}-\frac{11x^8}{8}+\frac{11x^6}{6}
-\frac{11x^4}{8}+\frac{5x^2}{12},\\
\int x^{12}&=&\frac{x^{13}}{13}+\frac{x^{12}}{2}+x^{11}-\frac{11x^9}{6}+\frac{22x^7}{7}-\frac{33x^5}{10}+\frac{5x^3}{3}-\frac{691x}{2730},\\
\int x^{13}&=&\frac{x^{14}}{14}+\frac{x^{13}}{2}+\frac{13x^{12}}{12}
-\frac{143x^{10}}{60}+\frac{143x^8}{28}-\frac{143x^6}{20}+\frac{65x^4}{12}
-\frac{691x^2}{420},\\
\int x^{14}&=&\frac{x^{15}}{15}+\frac{x^{14}}{2}+\frac{7x^{13}}{6}-\frac{91x^{11}}{30}
+\frac{143x^9}{18}-\frac{143x^7}{10}+\frac{91x^5}{6}-\frac{691x^3}{90}
+\frac{7x}{6},\\
\int x^{15}&=&\frac{x^{16}}{16}+\frac{x^{15}}{2}+\frac{5x^{14}}{4}
-\frac{91x^{12}}{24}+\frac{143x^{10}}{12}-\frac{429x^8}{16}
+\frac{455x^6}{12}-\frac{691x^4}{24}+\frac{35x^2}{4},\\
\int x^{16}&=&\frac{x^{17}}{17}+\frac{x^{16}}{2}+\frac{4x^{15}}{3}
-\frac{14x^{13}}{3}+\frac{52x^{11}}{3}-\frac{143x^9}{3}
+\frac{260x^7}{3}-\frac{1382x^5}{15}\\
&&+\frac{140x^3}{3}-\frac{3617x}{510}.
\end{eqnarray*}

24. But if on the other hand $x$ does not always have positive exponents in the
general term of a series, then too the expression will come out to be a sum
of infinitely many terms, because series of this kind do not admit general
summation, but rather involve quadratures. Still though, I have observed
that by means of this formula this type of series can easily be summed
very closely. This has a great utility in series which converge slowly
and others that are difficult to sum. I will explain through examples how this
is done.

25. Thus before any others I will first consider the harmonic series
\[
1+\frac{1}{2}+\frac{1}{3}+\frac{1}{4}+\textrm{etc.},
\]
whose general term is $\frac{1}{x}$, and let $S$ be the summatory term which is sought. So it is $X=\frac{1}{x}$ and
\[
\int Xdx=\textrm{Const.}+lx.
\]
And then
\[
\frac{dX}{dx}=\frac{-1}{x^2}, \quad \frac{d^3X}{dx^3}=\frac{-1\cdot 2\cdot 3}{x^4}, \quad \frac{d^5X}{dx^5}=\frac{-1\cdot 2\cdot 3\cdot 4\cdot 5}{x^6}
\quad \textrm{etc.}
\]
Substituting these yields
\begin{eqnarray*}
S&=&\textrm{Const.}+lx+\frac{1}{2x}-\frac{1}{12x^2}+\frac{1}{120x^4}
-\frac{1}{252x^6}+\frac{1}{240x^8}-\frac{1}{132x^{10}}\\
&&+\frac{691}{32760x^{12}}-\frac{1}{12x^{14}}+\textrm{etc.},
\end{eqnarray*}
where the constant that is added needs to be such that by putting $x=0$ it makes
$S=0$.
Certainly however the constant cannot be determined from this, because all
the terms are infinitely large.

26. Indeed to determine the constant another case should be considered, in
which the sum of the series is known; this can be obtained if a certain
number of terms are gathered into a single sum. Therefore let us add
the first $10$ terms
\[
1+\frac{1}{2}+\frac{1}{3}+\ldots+\frac{1}{10}
\]
whose sum turns out to be
\[
=2,9289682539682539;
\]
this should be equal to the sum of the terms from the formula, namely
\[
\textrm{Const.}+l10+\frac{1}{20}-\frac{1}{1200}+\frac{1}{1200000}-\frac{1}{252000000}
+\frac{1}{24000000000}
-\frac{1}{1320000000000}+\textrm{etc.}
\]
With this done, because one finds that
\[
l10=2,302585092994045684
\]
the added constant will be
\[
=0,5772156649015329
\]
and with this determined once, any sum of terms of this series can be found.

27. I have investigated the sum of $100,1000,10000$ etc. terms of the series
$1+\frac{1}{2}+\frac{1}{3}+\frac{1}{4}+\frac{1}{5}+$ etc. by this rule,
and I have found the following, 
\begin{eqnarray*}
\int 10&=&2,9289682539682539,\\
\int 100&=&5,1873775176396203,\\
\int 1000&=&7,4854708605503449,\\
\int 10000&=&9,7876060360443823,\\
\int 100000&=&12,0901461298634280,\\
\int 1000000&=&14,3927267228657236.
\end{eqnarray*}

28. If the first term of the series, 1, is taken, it will be $S=1$ and $x=1$,
and hence $lx=0$. From the equation we therefore get
\[
0,4227843350984670=\frac{1}{2}-\frac{1}{12}+\frac{1}{120}-\frac{1}{252}
+\frac{1}{240}-\frac{1}{132}+\frac{691}{32760}-\frac{1}{12}+\textrm{etc.}
\]
This series is very irregular and not even convergent, and the sum is
found only approximately. However the sum of the series continued to
infinity will be
\[
=l\infty+0,5772156649015329,
\]
which happens by putting $x=\infty$.

29. Now let us proceed to considering this series
\[
1+\frac{1}{3}+\frac{1}{5}+\frac{1}{7}+\frac{1}{9}+\textrm{etc.}
\]
in which $X=\frac{1}{2x-1}$ and
\[
\int Xdx=\textrm{Const.}+\frac{1}{2}l(2x-1)
\]
and also
\[
\frac{dX}{dx}=\frac{-2}{(2x-1)^2}, \quad \frac{d^3X}{dx^3}=\frac{-2\cdot 4\cdot 6}{(2x-1)^4},\quad \frac{d^5X}{dx^5}=\frac{-2\cdot 4\cdot 6\cdot 8\cdot 10}{(2x-1)^6} \quad \textrm{etc.}
\]
With these found, the sum of the proposed series will be
\begin{eqnarray*}
S&=&\textrm{Const.}+\frac{1}{2}l(2x-1)+\frac{1}{2(2x-1)}-\frac{1}{6(2x-1)^2}
+\frac{1}{15(2x-1)^4}-\frac{8}{63(2x-1)^6}\\
&&+\frac{8}{15(2x-1)^8}-\frac{128}{32(2x-1)^{10}}+\frac{256\cdot 691}{4095(2x-1)^{12}}
-\frac{2048}{3(2x-1)^{14}}+\frac{1024\cdot 3617}{255(2x-1)^{16}}-\textrm{etc.}
\end{eqnarray*}

30. The constant quantity in this case cannot be determined as easily
as that in the previous case by actual addition of several terms.
Indeed in this case a great help is that this constant can be determined
from the preceding one. Namely the sum of the series
\[
1+\frac{1}{3}+\frac{1}{5}+\frac{1}{7}+\textrm{etc.}
\]
continued to infinity is $=\textrm{Const.}+\frac{1}{2}l\infty$.
Let us subtract the harmonic series from twice this series; we will have
\[
1-\frac{1}{2}+\frac{1}{3}-\frac{1}{4}+\textrm{etc.},
\]
whose sum it turns out is $l2$. Therefore it will be
\[
l2=2\textrm{const.}+l\infty-l\infty-0,577215\, \textrm{etc.}
\]
and hence the constant that is sought is
\[
=0,6351814227307392.
\]

31. I proceed now to more complicated series, and I consider
\[
1+\frac{1}{4}+\frac{1}{9}+\frac{1}{16}+\textrm{etc.}
\]
the reciprocals of the squares, whose general term is $\frac{1}{x^2}=X$.
Therefore it will be
\[
\int Xdx=\textrm{Const.}-\frac{1}{x}
\]
and
\[
\frac{dX}{dx}=\frac{-2}{x^3},\quad \frac{d^3X}{dx^3}=\frac{-2\cdot 3\cdot 4}{x^5},
\quad \frac{d^5X}{dx^5}=\frac{-2\cdot 3\cdot 4\cdot 5\cdot 6}{x^7}\quad \textrm{etc.}
\]
With these substituted it will be
\[
1+\frac{1}{4}+\frac{1}{9}+\frac{1}{16}+\ldots+\frac{1}{x^2}=S
\]
\[
=\textrm{Const.}-\frac{1}{x}+\frac{1}{2x^2}-\frac{1}{6x^3}+\frac{1}{30x^5}
-\frac{1}{42x^7}+\frac{1}{30x^9}-\frac{5}{66x^{11}}+\frac{691}{2730x^{13}}
-\frac{7}{6x^{15}}+\textrm{etc.}
\]
where the constant quantity should be determined from a special case.

32. Thus I actually added the first ten terms of this series and I found
that their sum is 
\[
1,549767731166540.
\]
Since in this case $x=10$, if this is added to
\[
\begin{split}
&\frac{1}{10}-\frac{1}{200}+\frac{1}{6000}-\frac{1}{3000000}+\frac{1}{420000000}
-\frac{1}{30000000000}+\frac{1}{1320000000000}\\
&-\frac{691}{27300000000000000}
+\frac{7}{6000000000000000}-\textrm{etc.},
\end{split}
\]
one gets an added constant $=1,64493406684822643647$.
And this constant is equal to the sum of the series continued to infinity;
for by putting $x=\infty$ it will be $S=\textrm{Const.}$, with all the terms
vanishing.

33. In a similar way for the reciprocals of the cubes
\[
1+\frac{1}{8}+\frac{1}{27}+\frac{1}{64}+\textrm{etc.}
\]
if the first ten terms are added, this sum is obtained
\[
1,197531985674193.
\]
Whence one finds that the constant which should be added in the summation
of this series is
\[
=1,202056903159594.
\]
And this number is equal to the sum of the series
\[
1+\frac{1}{8}+\frac{1}{27}+\frac{1}{64}
\]
continued to infinity.

And for the biquadrates
\[
1+\frac{1}{16}+\frac{1}{81}+\textrm{etc.}
\]
the sum is
\[
=1,0823232337110824.
\]

34. Let us now consider by this method the series by which the area of
the circle whose diameter is $1$ is exhibited, namely
\[
1-\frac{1}{3}+\frac{1}{5}-\frac{1}{7}+\frac{1}{9}-\textrm{etc.}
\]
or
\[
\frac{2}{1\cdot 3}+\frac{2}{5\cdot 7}+\frac{2}{9\cdot 11}+\frac{2}{13\cdot 15}+\textrm{etc.},
\]
whose general term is
\[
\frac{2}{(4x-3)(4x-1)}
\]
or resolving into factors
\[
\frac{1}{4x-3}-\frac{1}{4x-1}.
\]
For finding the approximate sum of this series,
\[
X=\frac{1}{4x-3}-\frac{1}{4x-1}
\]
and
\[
\int Xdx=\textrm{Const.}-\frac{1}{4}l\frac{4x-1}{4x-3}
\]
and then
\[
\frac{dX}{dx}=\frac{-4}{(4x-3)^2}+\frac{4}{(4x-1)^2},
\quad \frac{d^3X}{dx^3}=\frac{-4\cdot 8\cdot 12}{(4x-3)^4}+\frac{4\cdot 8\cdot 12}{(4x-1)^4} \quad \textrm{etc.}
\]
From this the sum of the series
\[
\frac{2}{1\cdot 3}+\frac{2}{5\cdot 7}+\ldots+\frac{2}{(4x-3)(4x-1)}+\textrm{etc.}
\]
will be
\begin{eqnarray*}
S&=&\textrm{Const.}-\frac{1}{4}l\frac{4x-1}{4x-3}+\frac{1}{2}\bigg(\frac{1}{4x-3}
-\frac{1}{4x-1}\bigg)-\frac{1}{3}\bigg(\frac{1}{(4x-3)^2}-\frac{1}{(4x-1)^2}\bigg)\\
&&+\frac{8}{15}\bigg(\frac{1}{(4x-3)^4}-\frac{1}{(4x-1)^4}\bigg)
-\frac{256}{63}\bigg(\frac{1}{(4x-3)^6}-\frac{1}{(4x-1)^6}\bigg)\\
&&+\frac{1024}{15}\bigg(\frac{1}{(4x-3)^8}-\frac{1}{(4x-1)^8}\bigg)
-\frac{4^8}{33}\bigg(\frac{1}{(4x-3)^{10}}-\frac{1}{(4x-1)^{10}}\bigg)+\textrm{etc.}
\end{eqnarray*}
Truly even if ten terms of this series are added it will not converge enough
so that a proper constant could be exhibited. 
But four times the constant is equal to the periphery of a circular
whose diameter is $=1$.

\end{document}